# An Improved Benders Decomposition Algorithm for Steady-State Dispatch Problem in an Integrated Electricity-Gas System

Han Gao, *Student Member, IEEE*, Zhengshuo Li, *Member, IEEE*

*Abstract*—Optimally operating an integrated electricity-gas system (IEGS) is significant for the energy sector. However, the IEGS operation model's nonconvexity makes it challenging to solve the optimal dispatch problem in the IEGS. This letter proposes an improved Benders decomposition (IBD) algorithm catering to a commonly used steady-state dispatch model of the IEGS. This IBD algorithm leverages a refined decomposition structure where the subproblems become linear and ready to be solved in parallel. We analytically compare our IBD algorithm with an existing Benders decomposition algorithm and a typical piecewise linearization method. Case studies have substantiated the higher computational efficiency of our IBD algorithm.

## I. Introduction

WITH the development of gas-fired generators and power to gas (P2G) facilities, the interdependency between power systems and natural gas systems has increased [1]. Optimally operating an integrated electricity-gas system (IEGS) provides a promising way to enhance overall energy utilization, and this positive effect has been confirmed in the literature, e.g., [1]-[3]. However, due to the nonconvexity, e.g., in Weymouth gas flow equations, solving the optimal dispatch model of an IEGS is computationally challenging. Piecewise linearization (PWL) methods entailing auxiliary binary variables are widely employed [4], but the accuracy depends on the number of linearization segments. To achieve a satisfying solution, one needs to introduce abundant auxiliary binary variables and constraints, resulting in an expensive computational cost. Alternatively, [3] and [5] present a Benders decomposition (BD) algorithm where the steady-state natural gas system constraints and the power system dispatch problem are handled in a master-subproblem framework. This letter proposes an improved Benders decomposition (IBD) algorithm that caters to a commonly used IEGS dispatch model. Case studies demonstrate that our IBD algorithm has a significantly higher computational efficiency than the preceding BD algorithm as well as regular PWL methods.

## II. Steady-State Dispatch Model of an IEGS

Below is the model of a steady-state IEGS dispatch problem that is commonly adopted in the literature.

This work was supported by National Natural Science Foundation of China under Grant 52007105. Han Gao and Zhengshuo Li are with the School of Electrical Engineering, Shandong University, Jinan, China. Zhengshuo Li is the corresponding author (email: zsli@sdu.edu.cn).

### A. Steady-state Power System Dispatch Model

A standard power system dispatch model is formulated as

$$-\boldsymbol{P}_{tra} \leq \mathbf{H}_g \boldsymbol{P}_{g,t} - \mathbf{H}_d \boldsymbol{P}_{d,t} \leq \boldsymbol{P}_{tra}, \tag{1}$$

$$\sum_{i\in\Omega_F} P_{f,i}^t + \sum_{j\in\Omega_G} P_{g,j}^t + \sum_{n\in\Omega_W} P_{w,n}^t = \sum_{m\in\Omega_D} P_{d,m}^t + \sum_{k\in\Omega_P} P_{p,k}^t, \tag{2}$$

$$P_{f,i}^{\min} \leq P_{f,i}^t \leq P_{f,i}^{\max}, i\in\Omega_F, P_{g,j}^{\min} \leq P_{g,j}^t \leq P_{g,j}^{\max}, j\in\Omega_G, \tag{3}$$

$$0 \leq P_{w,n}^t \leq P_{w,n}^{t,\max}, n\in\Omega_W, \tag{4}$$

$$R_{f,i}^{dw} \leq P_{f,i}^t - P_{f,i}^{t-1} \leq R_{f,i}^{up}, R_{g,j}^{dw} \leq P_{g,j}^t - P_{g,j}^{t-1} \leq R_{g,j}^{up}, \tag{5}$$

$$F_{f,i}^t = a_{f,i}(P_{f,i}^t)^2 + b_{f,i} P_{f,i}^t + c_{f,i}, i\in\Omega_F, \tag{6}$$

$$G_{g,j}^t = a_{g,j}(P_{g,j}^t)^2 + b_{g,j} P_{g,j}^t + c_{g,j}, j\in\Omega_G, \tag{7}$$

where $\Omega_F$, $\Omega_G$, $\Omega_W$, $\Omega_D$, $\Omega_P$ denote the sets of coal-fired units, gas-fired units, wind farms, electricity loads, and P2G facilities, respectively; $P_{f,i}^t$, $P_{g,j}^t$, and $P_{w,n}^t$ are the scheduled outputs of coal-fired unit $i$, gas-fired unit $j$, and wind farm $n$ at time $t$, respectively; $P_{d,m}^t$, $P_{p,k}^t$ denote the consumed electricity of load $m$ and P2G facility $k$ at time $t$; $P_{f,i}^{\min}$, $P_{f,i}^{\max}$, $P_{g,j}^{\min}$, $P_{g,j}^{\max}$ are the lower and upper generation limits of coal-fired unit $i$ and gas-fired unit $j$, respectively; $R_{f,i}^{dw}$, $R_{f,i}^{up}$, $R_{g,j}^{dw}$, $R_{g,j}^{up}$ are the downward and upward ramping limits of coal-fired unit $i$ and gas-fired unit $j$, respectively; $P_{w,n}^{t,\max}$ is the available power of wind farm $n$ at time $t$; $F_{f,i}^t$ and $G_{g,j}^t$ are the cost of coal-unit $i$ and gas consumption of gas-fired unit $j$, respectively, and $a_{f,i}$, $b_{f,i}$, $c_{f,i}$, $a_{g,j}$, $b_{g,j}$, $c_{g,j}$ are the coefficients; vectors $\boldsymbol{P}_{g,t}$ and $\boldsymbol{P}_{d,t}$ denote the output of generation units and load across the power system at time $t$, respectively; $\boldsymbol{P}_{tra}$ is the transmission capacity of all lines; $\mathbf{H}_g$, $\mathbf{H}_d$ are the matrices of the distribution factors of line flow regarding the generation and load buses, respectively.

The constraint (1) denotes the transmission capacity limits of all lines; (2) denotes the system power balance constraint; (3)-(4) limit the power output of the units; (5) denotes the ramping rate limits of the coal-fired and gas-fired units; (6), (7) describe the consumption functions of coal-fired unit $i$ and gas-fired unit $j$, respectively.

### B. Steady-state Natural Gas System Model

$$G_{ij}^t \left| G_{ij}^t \right| = k_l \left( \pi_{i,s}^t - \pi_{j,s}^t \right),\ l\in\Omega_L, i,j\in\Omega_{Nl}, \tag{8}$$

$$\pi_{i,s}^{\min} \leq \pi_{i,s}^t \leq \pi_{i,s}^{\max},\ i\in\Omega_N, \tag{9}$$

$$r_m^{\min} \pi_{i,s}^t \leq \pi_{j,s}^t \leq r_m^{\max} \pi_{i,s}^t, m\in\Omega_C, \tag{10}$$

$$G_{c,m}^t = \alpha_m^t G_m^t, m\in\Omega_C, \tag{11}$$

$$G_{s,i}^{\min} \leq G_{s,i}^t \leq G_{s,i}^{\max}, i\in\Omega_S, -G_z^{out} \leq G_{r,z}^t \leq G_z^{in}, z\in\Omega_R, \tag{12}$$

$$S_z^t = S_z^{t-1} + G_{r,z}^t, S_z^{\min} \leq S_z^t \leq S_z^{\max}, S_z^t = S_z^{t+24}, t\in T, \tag{13}$$

$$G_{p,k}^t = \eta_k P_{p,k}^t, 0 \leq P_{p,k}^t \leq P_{p,k}^{\max}, k\in\Omega_P, \tag{14}$$

$$\sum_{i\in\Omega_{Sn}} G_{s,i}^t + \sum_{k\in\Omega_{Pn}} G_{p,k}^t - \sum_{m\in\Omega_{cn}} G_{c,m}^t - \sum_{j\in\Omega_{Nn}} G_{nj}^t - \sum_{z\in\Omega_{Rn}} G_{r,z}^t$$
$$+ \sum_{m\in\Omega_{Cn}} s_{n,m} G_m^t = \sum_{j\in\Omega_{Gn}} G_{g,j}^t + \sum_{w\in\Omega_{Dn}} G_{d,w}^t, \quad \forall n \in \Omega_N, \quad (15)$$

where $\Omega_L$, $\Omega_N$, $\Omega_C$, $\Omega_S$, $\Omega_R$ denote the sets of the pipelines without compressors, gas system nodes, compressors, gas wells, gas storage devices, respectively; $\Omega_{Nl}$ includes the two ends of the pipeline $l$; $\Omega_{Sn}$, $\Omega_{Pn}$, $\Omega_{Rn}$, $\Omega_{Cn}$, $\Omega_{Gn}$ denote the subsets of $\Omega_S$, $\Omega_P$, $\Omega_R$, $\Omega_C$, $\Omega_G$ where the elements are connected to node $n$; $\Omega_{Dn}$ denotes the set of the gas load at node $n$; $\Omega_{cn}$ is the set of compressors which transmit gas from node $n$; $\Omega_{Nn}$ is the set of nodes connected to node $n$ by a pipeline without compressors; $T$ denotes the entire scheduling time horizon with $N_T$ periods; $G_{ij}^t$ is the gas flow through the pipeline from node $i$ to node $j$ at time $t$; $k_l$ is the gas constant of pipeline $l$; $\pi_{i,s}^t$ is the squared pressure of node $i$ at time $t$; $\pi_{i,s}^{\min}$, $\pi_{i,s}^{\max}$ are the lower and upper squared pressure limits on node $i$, respectively; $r_m^{\min}$, $r_m^{\max}$ denote the minimum and maximum squared compression ratio of compressor $m$, respectively; $G_m^t$, $\alpha_m^t$, and $G_{c,m}^t$ are the gas flow through compressor $m$, consumption coefficient, and gas consumption of compressor $m$ at time $t$, respectively; $G_{s,i}^{\min}$, $G_{s,i}^{\max}$, $G_{s,i}^t$ are the lower and upper output limits and the output of gas well $i$, respectively; $G_z^{out}$, $G_z^{in}$, $S_z^{min}$, $S_z^{max}$, and $S_z^t$ are the maximum extraction and injection rate, lower and upper gas storage level, and storage level of storage device $z$, respectively; $G_{r,z}^t$ is the gas exchanging rate of storage device $z$ at time $t$: it is positive for injection and negative for extraction; $\eta_k$, $P_{p,k}^{\max}$ are the conversion coefficient and upper limit of the electricity demand of P2G facility $k$, respectively; $G_{p,k}^t$, $G_{d,w}^t$ denote the gas outputs of P2G facility $k$ and gas load $w$ at time $t$, respectively; $s_{n,m}$ is the direction indicator: it is 1 if $n$ is the outlet node of compressor $m$, and -1 otherwise.

As for a pipeline without compressors, the Weymouth equation (8) shows the relationship between the nodal pressure and the transmitted gas flow. As for a pipeline with a compressor, the gas transmitted from the sending end is $G_m^t+G_{c,m}^t$, while only $G_m^t$, which is usually deemed free from the constraint (8) due to the compressor's functioning [5], is transmitted to the receiving end. $G_{c,m}^t$ is approximately proportional to $G_m^t$ in (11), and $\alpha_m^t$ depending on several parameters can be further simplified as a constant considering that a compressor's gas consumption ratio is typically 3%-5% [2]. Next, (9) and (10) limit the range of the nodal pressure and the compression ratios; (12) limits the output of gas well $i$ and the rate of gas injection and extraction from storage $z$; (13) describes the gas storage level limit of gas storage $z$ at time $t$. As per [2], gas storage is typically considered high-efficiency with negligible losses during the gas injection and extraction process, justifying using a single variable $G_{r,z}^t$ to stand for the gas exchanging rate for storage $z$ at time $t$; (14) denotes the energy conversion efficiency of the P2G facilities [6]; (15) denotes the gas flow balance of node $n$ of the natural gas system.

### C. Objective of IEGS

The objective is to minimize the operating costs of the IEGS as well as the wind curtailment penalty, formulated as

$$\min f = \sum_{t\in T}\left(\sum_{i\in\Omega_F} F_{f,i}^t + \sum_{j\in\Omega_S} b_{s,j} G_{s,j}^t + \sum_{m\in\Omega_W} F_{w,m}^t\right), \quad (16)$$

$$F_{w,m}^t = \rho\left(P_{w,m}^{t,\max} - P_{w,m}^t\right), m \in \Omega_W, \quad (17)$$

where $b_{s,j}$ denotes the cost coefficient of gas well $j$; $F_{w,m}^t$ is the wind curtailment cost; $\rho$ is the penalty factor and set to 35\$/MW [7] in the case study.

### III. IMPROVED BENDERS DECOMPOSITION ALGORITHM

The above model is nonconvex due to (6)-(8). However, it is easy to see and also a regular practice that (6) can be relaxed as a convex inequality, which is tight at the optimal solution both because $a_{f,i}$ is positive and $F_{f,i}^t$ only appears in the objective function. Since $G_{g,j}^t$ in (7) also appears in (15), the preceding relaxation technique, if applied to (7), could lose its tightness, so we resort to the PWL method in [8] to linearize (7). However, further applying the PWL method to (8) will introduce so many additional binaries to incur heavy computational burdens, as shown in the following case study.

Alternatively, with (6), (7) handled as alluded to earlier but (8) remaining as it is, the BD algorithm in [3], [5] can be employed to solve this reformulated dispatch model. However, it is easy to see that this BD algorithm's subproblem is nonconvex with fixed master-problem variables, violating the assumptions about the classical generalized BD algorithm [9]. Consequently, it may converge slowly or even fail to converge.

To overcome this issue, we design an IBD algorithm of a refined decomposition structure that conforms to the assumptions in [9]. From a computational standpoint, different decomposition structures in BD algorithms may directly impact the computational performance [10]. Hence, it is reasonable to expect our IBD would perform better than the preceding BD algorithm as it conforms to the assumption in [9]. Moreover, as will be seen, our IBD's subproblem is ready to be solved in parallel, further accelerating the computation. The IBD's master and subproblems are briefly shown below, and the algorithm flowchart is shown in Fig.1(a).

- **Master Problem**: It encompasses the (reformulated) (1)-(7), (11)-(17) and the accumulated cuts generated by the subproblem in the iteration.
- **Gas Flow Feasibility Check Subproblem**: This subproblem checks the gas flow feasibility with fixed master-problem variables. Its formulation for each time $t$ is shown as

$$\min_{\pi_{is}^t, L1_{tl}, L2_{tl}} g^t = \sum_{l\in\Omega_L}(L1_{tl} + L2_{tl}), \quad (18)$$

subject to (9)-(10) and $\overline{G}_{ij}^t\left|\overline{G}_{ij}^t\right| - k_l(\pi_{is}^t - \pi_{js}^t) = L1_{tl} - L2_{tl}, l\in\Omega_L, \quad (19)$

where $\overline{G}_{ij}^t$ denotes the gas flow specified by the master problem; $L1_{tl}, L2_{tl}$ are the nonnegative slack variables. The objective (18) is to minimize the sum of possible mismatches in the Weymouth gas flow equation with a specified $\overline{G}_{ij}^t$.

Let $u_l^t$ denote the dual variables related to (19). One can generate the infeasibility cut [5] below for each time $t$:

$$g^t + \sum_{l\in\Omega_L}(\sigma_l^t u_l)(G_{ij}^t - \overline{G}_{ij}^t) \leq 0, \quad (20)$$

where $\sigma_l^t$ equals $2\left|\overline{G}_{ij}^t\right|$, i.e., the gradient of $G_{ij}^t\left|G_{ij}^t\right|$ where $G_{ij}^t$ is equal to $\overline{G}_{ij}^t$.

The structural differences between the preceding BD and our IBD algorithm are summarized below. First, the BD algorithm has both infeasibility and optimality cuts [3]. However, as for the IBD, since the IEGS's objective function is holisti-cally optimized in the master problem, only the infeasibility cuts are generated whenever positive minimal mismatches appear. Second, the subproblem of the BD algorithm is nonconvex while that of the IBD is linear given a fixed $\overline{G}_{ij}^t$. Third, the IBD's subproblem is decomposable in terms of the timestamp $t$ to be solved in parallel and then generate multiple cuts. This point is quantitatively shown in Table I, where $N_c$, $N_n$, $N_r$, $N_l$, and $N_s$ are the number of the compressors, gas nodes, storage devices, pipelines without compressors, and gas wells, respectively.

TABLE I. COMPARISON BETWEEN SUBPROBLEMS OF BD AND IBD

| | Property | Number of Variables[a] | Number of Cuts |
|---|---|---|---|
| BD | Nonconvex | $N_T(2N_c+N_n+2N_r+N_l+N_s)$ | 1 |
| IBD | Linear | $N_T$ subproblems; number of each subproblem variables is $(2N_l+N_n)$ | $N_T$ |

## IV. CASE STUDIES

To compare the efficacy of different algorithms for a multi-period dispatch problem, we considered a 24-hour IEGS dispatch problem with the schedule time interval being one hour (so $N_T$ is 24). By modifying the system in [1], three IEGSs were constructed by coupling a modified IEEE 24-bus power system with three natural gas systems of different scales, respectively. The P2G facility is installed at bus 8 of the power system, and its maximal power is 200MW. All the three cases below were tested on a PC with Intel Core i7 2.90GHz CPU and 16GB RAM. The codes were implemented in the MATLAB environment using YALMIP and GUROBI.

- CASE1: IEEE 24-bus power system is connected with an 8-pipeline natural gas system. Compared with [1], nodes 8, 9, 10, 11, 13, 14 of the gas system and their adjoining pipelines were removed; meanwhile, node 6 and 12 were connected by a pipeline.
- CASE2: Similar to CASE1 except that nodes 8, 9 and their adjoining pipelines were recovered.
- CASE3: Similar to CASE2 except that nodes 10, 11 and their adjoining pipelines were recovered further.

For each case, our IBD algorithm is compared with the PWL method in [8] and the BD algorithm in [3]. The results are shown in Table II. The iteration numbers of the BD and IBD are shown in Fig.1(b).

TABLE II. RESULTS OF PWL, BD, AND IBD IN THREE CASES

| Method | Solution time (sec.) | | | Operating cost (USD) | | |
|---|---|---|---|---|---|---|
| | PWL | BD | IBD | PWL | BD | IBD |
| CASE1 | 85.75 | 72.99 | 6.17 | 4183498 | 4216337 | 4183488 |
| CASE2 | 3291.64 | 118.84 | 9.79 | 4039586 | 4072246 | 4039582 |
| CASE3 | >7200 | 325.99 | 11.40 | -- | 4453075 | 4420736 |

Table II shows that the PWL consumes the longest time for all the cases because the PWL linearizes Weymouth equations for acceptable accuracy by using numerous integer variables and constraints. For example, in CASE1, one Weymouth equation is approximated by 56 linearization segments, leading to 8000+ integer variables and 8000+ continuous variables for this 24-hour dispatch problem. With an increasing number of pipelines, the number of integer variables increases, incurring the notorious curse of dimensionality.

Next, it can be seen that the IBD is faster than the BD and also yields lower operating costs. In principle, this is owing to its refined decomposition structure. To be more specific, two factors are analyzed below. First, the decomposed subproblems of the IBD generate $N_T$ cuts that are simultaneously added to the master problem in each iteration, but the BD generates only a single cut in each iteration—as per [10], adding multiple cuts is usually better than adding a single (even high-density) cut in reducing the number of iterations.

Second, the linear subproblems of the IBD can be efficiently solved by off-the-shelf solvers in parallel. Contrarily, the subproblem of the BD is nonconvex, calling for, e.g., a tailored Newton–Raphson algorithm [3]. However, when it is hard to specify a valid initial trial value, the Newton–Raphson algorithm might converge slowly or even diverge, as was observed in our study. Moreover, since the solution to the nonconvex subproblem of the BD is not surely globally optimal, the generated cuts might have improperly cut off the feasible region of the master problem, which would in turn call for more iterations and even lead to a solution inferior to those of the PWL and IBD, as is revealed in Table II and Fig.1(b).

To conclude, our IBD algorithm has a significantly higher computational efficiency than the BD and PWL methods in solving a steady-state IEGS dispatch problem, especially when the natural gas system has multiple pipelines.

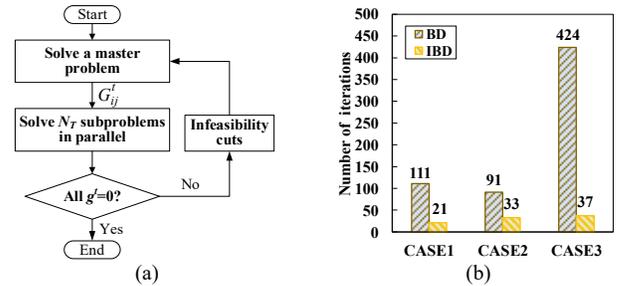

Fig. 1. (a) Algorithm flowchart of IBD. (b) Iteration numbers of BD and IBD.

---

[a] Here, for both BD and IBD algorithms, we exclude from the subproblem variables the gas-related variables that have been fixed in the master problem.